\documentclass[11pt,a4paper]{article}
\usepackage[T1]{fontenc}
\usepackage{lmodern,amsmath,amssymb,amsthm,mathtools,microtype,booktabs,array}
\usepackage[margin=25mm]{geometry}
\usepackage{graphicx,subcaption,placeins}
\usepackage[colorlinks=true,linkcolor=blue,citecolor=blue,urlcolor=blue]{hyperref}
\numberwithin{equation}{section}
\newtheorem{theorem}{Theorem}[section]

\theoremstyle{definition}
\theoremstyle{remark}\newtheorem{remark}[theorem]{Remark}
\newcommand{\ii}{\mathrm i}
\newcommand{\dd}{\,\mathrm d}
\newcommand{\R}{\mathbb R}
\newcommand{\cv}{c}
\newcommand{\cO}{\mathcal O}
\newcommand{\cI}{\mathcal I}
\newcommand{\cE}{\mathcal E}

\DeclareMathOperator{\arcosh}{arcosh}
\DeclareMathOperator{\Arg}{Arg}
\allowdisplaybreaks[2]
\title{Limit shape for domain-wall six-vertex model}

\makeatletter
\newcommand{\address}[2][]{\g@addto@macro\@author{\\[0.5em]\normalsize #2}}
\newcommand{\email}[1]{\g@addto@macro\@author{\\\normalsize\href{mailto:#1}{\texttt{#1}}}}
\makeatother

\author{Alexey Bufetov}
\date{}

\begin{document}
	\maketitle
	
	\begin{abstract}
		In this note we announce limit shape theorems for the height function of the six-vertex model in the antiferroelectric (for arbitrary $\Delta <-1$) and disordered (in a restricted range $-1 < \Delta \le 1/2$) regimes under domain-wall boundary conditions.  
	\end{abstract}
	
	\section{The model}\label{sec:model}
	
	\subsection{Vertices, weights, and boundary conditions}
	Let $N\ge1$. Place the six-vertex vertices at the square grid
	\[
	(i-\tfrac12,j-\tfrac12),\qquad 1\le i,j\le N,
	\]
	and join nearest neighbors, adding the $4N$ boundary half-edges to the sides of $[0,N]^2$. Each edge is either occupied or empty, which we label as $\{1,0\}$. At a vertex denote the west, north, east and south occupations by $(w,n;e,s)$. The local constraint of the six-vertex model is $w+n=e+s$. 
	
	Occupied edges are directed east horizontally and south vertically; thus $w,n$ are incoming and $e,s$ outgoing occupations. The six allowed states have weights
	\begin{equation}\label{eq:weights-table}
		\begin{array}{c|c|c}
			(w,n)&(e,s)&\text{weight}\\\hline
			(0,0)&(0,0)&a\\
			(1,1)&(1,1)&a\\
			(1,0)&(1,0)&b\\
			(0,1)&(0,1)&b\\
			(1,0)&(0,1)&\cv\\
			(0,1)&(1,0)&\cv
		\end{array}
		\qquad a,b,\cv>0.
	\end{equation} 
	

	The domain-wall boundary conditions are
	\begin{equation}\label{eq:DW-bits}
		\begin{aligned}
			w&=0\quad\text{on the west side},& n&=1\quad\text{on the north side},\\
			e&=1\quad\text{on the east side},& s&=0\quad\text{on the south side}.
		\end{aligned}
	\end{equation}
	Let $\Omega_N$ be the finite set of configurations satisfying these conditions. Writing $W_v(\omega)$ for the weight at vertex $v$, define
	\begin{equation}\label{eq:measure}
		Z_N=\sum_{\omega\in\Omega_N}\prod_v W_v(\omega),\qquad
		\mathbb P_N^{a,b,\cv}(\omega)=Z_N^{-1}\prod_vW_v(\omega).
	\end{equation}
	We denote the measure by $\mathbb P_N$.
	
	\subsection{Height function}
	The faces, including the boundary sectors, are indexed by (dual) vertices $(i,j)\in\{0,\ldots,N\}^2$. Set $H_N(0,0)=0$. Around the vertex $(i-\frac12,j-\frac12)$, use the adjacent face labels $\mathrm{SW},\mathrm{SE},\mathrm{NW},\mathrm{NE}$, and impose
	\begin{equation}\label{eq:height-local-rule}
		\begin{aligned}
			H_{\mathrm{SE}}-H_{\mathrm{SW}}&=s,&
			H_{\mathrm{NW}}-H_{\mathrm{SW}}&=w,\\
			H_{\mathrm{NE}}-H_{\mathrm{NW}}&=n,&
			H_{\mathrm{NE}}-H_{\mathrm{SE}}&=e.
		\end{aligned}
	\end{equation}
	It is easy to see that such a height function, with the fixed boundary value $H_N (0,0)=0$, is well defined and uniquely determines the configuration. 
	
	Under domain-wall boundary conditions, the boundary heights are
	\begin{equation}\label{eq:discrete-boundary-height}
		H_N(i,0)=H_N(0,j)=0,\qquad H_N(i,N)=i,\qquad H_N(N,j)=j.
	\end{equation}
	Set $h_N(i/N,j/N)=H_N(i,j)/N$. For each mesh square interpolate this function bilinearly. This determines a continuous rescaled height function $h_N (x,y)$ on $[0,1]^2$.
	
	\begin{figure}[tb]
		\centering
		\includegraphics[width=0.90\textwidth]{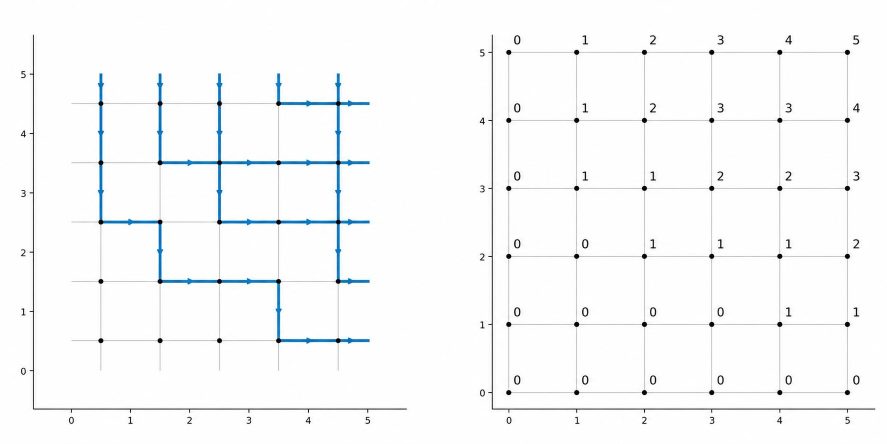}
		\caption{Example of a domain-wall configuration and its height function.}
		\label{fig:height-example}
	\end{figure}
	
	\subsection{Limit shape}
	A \textit{limit shape} is a deterministic $h_*: [0,1]^2 \to \R$ such that $\sup_{(x,y) \in [0,1]^2} | h_N (x,y)-h_*(x,y) | \to 0$ in probability (with respect to  $\mathbb P_N$), as $N \to \infty$. The four frozen affine pieces, labeled by their corners, are
	\begin{equation}\label{eq:frozen-planes}
		L_{\mathrm{SW}}=0,\qquad L_{\mathrm{SE}}=y,\qquad
		L_{\mathrm{NE}}=x+y-1,\qquad L_{\mathrm{NW}}=x.
	\end{equation}
	Theorems~\ref{thm:AF}-\ref{thm:D} and Section \ref{sec:curves} specify the open regions on which these pieces occur. They also specify the open liquid region and, in the antiferroelectric case, the open antiferroelectric region. An \emph{arctic curve} is a boundary between two of those regions.
	
	\subsection{Parameters}
	The model depends on two parameters, defined as
	\begin{equation}\label{eq:model-parameters}
		r=\frac ab>0,\qquad \Delta=\frac{a^2+b^2-\cv^2}{2ab}.
	\end{equation}
	For $\Delta<-1$ (AF regime) set
	\begin{equation}\label{eq:AFparams}
		\zeta=\arcosh(-\Delta),\qquad
		u=\frac12\log\frac{r+e^\zeta}{r+e^{-\zeta}},\qquad T=\zeta-u.
	\end{equation}
	Then $\zeta>0$, $0<u<\zeta$, and we have the following proportion:
	\begin{equation}\label{eq:AFweights}
		(a,b,\cv)\propto(\sinh T,\sinh u,\sinh\zeta),\qquad \Delta=-\cosh\zeta.
	\end{equation}
	For $-1<\Delta\le1/2$ (D regime) set
	\begin{equation}\label{eq:Dparams}
		\gamma=\arccos\Delta\in[\pi/3,\pi),\qquad \delta=\pi-\gamma,
	\end{equation}
	\begin{equation}\label{eq:D-u}
		u=\Arg(r-\cos\gamma+\ii\sin\gamma)\in(0,\delta),\qquad T=\delta-u.
	\end{equation}
	The weights satisfy the proportion 
	\begin{equation}\label{eq:Dweights}
		(a,b,\cv)\propto(\sin T,\sin u,\sin\delta),\qquad \Delta=\cos\gamma.
	\end{equation}
	We set $m=\zeta$ in the antiferroelectric case and $m=\delta$ in the disordered case. One has $0<u<m$.
	
	\section{Description of the answer}\label{sec:fredholm}
	
	\subsection{Determinant functional}
	All integral operators below act on the complex Hilbert space $L^2(I_B)$, where $I_B=[-B,B]$ with the Lebesgue measure. Define the signed kernel
	\begin{equation}\label{eq:signed-kernel}
		k(z):=
		\begin{cases}
			-\dfrac{\sinh(2\zeta)}{2\pi(\cosh(2\zeta)-\cos z)},&\text{AF case},\quad0<B<\pi,\\[6pt]
			\dfrac{\sin(2\gamma)}{2\pi(\cosh z-\cos(2\gamma))},&\text{D case},\quad B>0.
		\end{cases}
	\end{equation}
	In the second case set $k=0$ when $\gamma=\pi/2$. Define
	\begin{equation}\label{eq:base-operator}
		(\mathcal K_Bg)(s)=\int_{-B}^{B}k(s-t)g(t)\dd t,
		\qquad \mathfrak d_B=\det(I-\mathcal K_B).
	\end{equation}
	For $f\in L^2(I_B)$ define also the rank-one operator
	\begin{equation}\label{eq:rankone}
		(\mathcal U_fg)(s)=\int_{-B}^{B}f(t)g(t)\dd t\qquad(s\in I_B).
	\end{equation}
	Define the key determinant functional via
	\begin{equation}\label{eq:det-functional}
		\quad\mathcal D_B[f]=
		\frac{\det(I-\mathcal K_B+\mathcal U_f)}{\det(I-\mathcal K_B)}-1.
	\end{equation} 
	
	\subsection{Regime-specific functions}
	
	In the AF regime set 
	\begin{equation}\label{eq:AFkernel}
		A(z)=A_\zeta(z)=\frac{\sinh\zeta}{2\pi(\cosh\zeta-\cos z)},
	\end{equation}
	\begin{equation}\label{eq:AFell}
		\ell_v(z)=\frac{\zeta-v}{2\pi\zeta}
		+\frac1\pi\sum_{n\ge1}e^{-\zeta n}
		\frac{\sinh((\zeta-v)n)}{\sinh(\zeta n)}\cos(nz),\qquad0<v<\zeta.
	\end{equation}
	The series is holomorphic for $|\operatorname{Im}z|<\zeta+v$. In the D regime set
	\begin{equation}\label{eq:Dkernel}
		A(z)=A_\gamma(z)=\frac{\sin\gamma}{2\pi(\cosh z+\cos\gamma)},
	\end{equation}
	\begin{equation}\label{eq:Dell}
		\ell_v(z)=\frac{\sin(\delta-v)}{2\pi(\cosh z-\cos(\delta-v))}
		-\frac{\sin(\pi(\delta-v)/\delta)}
		{2\delta(\cosh(\pi z/\delta)-\cos(\pi(\delta-v)/\delta))},
		\quad0<v<\delta.
	\end{equation}
	Note that the function $\ell_v(z)$ is holomorphic for $|\operatorname{Im}z|<\delta+v$.
	
	For either regime, $0<v<m$, and $|\theta|<m-v$, define
	\begin{equation}\label{eq:local-determinants}
		\begin{aligned}
			p_v(B,\theta)&=\mathcal D_B\bigl[A(\,\cdot+\ii(\theta-v))\bigr],\\
			q_v(B,\theta)&=\mathcal D_B\bigl[A(\,\cdot+\ii(\theta+v))\bigr],\\
			R_v(B,\theta)&=\mathcal D_B\bigl[-\ell_v(\,\cdot+\ii\theta)\bigr].
		\end{aligned}
	\end{equation}
	All these three values are real.
	
	\subsection{Periodic extension}
	Let us set $\mathbb T_m=\R/(4m\mathbb Z)$, so functions of $\tau\in\mathbb T_m$ are $4m$-periodic. Define the parameter spaces via
	\begin{equation}\label{eq:parameter-spaces}
		\mathcal P_{\rm AF}=(0,\pi)\times\mathbb T_\zeta,
		\qquad \mathcal P_{\rm D}=(0,\infty)\times\mathbb T_\delta.
	\end{equation}
	Choose the representative $\tau\in[-T,3m+u)$. Away from the four joining points, define $P,Q,R$ by
	\begin{equation}\label{eq:charts}
		\begin{array}{c|c|c|c}
			\tau\text{ interval}&v&\theta&(P,Q,R)\\\hline
			(-T,T)&u&\tau&(p_u,q_u,R_u)\\
			(T,m+u)&T&m-\tau&(p_T,1-q_T,-p_T-R_T)\\
			(m+u,3m-u)&u&\tau-2m&(1-p_u,1-q_u,p_u+q_u+R_u-1)\\
			(3m-u,3m+u)&T&3m-\tau&(1-p_T,q_T,-q_T-R_T).
		\end{array}
	\end{equation}
	Each function in the last column is evaluated at $(B,\theta)$ from the corresponding parameter space. At a joining point of parameter intervals one uses the common one-sided limits.
	
	\subsection{Critical point}
	For $(x,y)\in[0,1]^2$ define
	\begin{equation}\label{eq:phase}
		\Phi_{x,y}(B,\tau)=xP(B,\tau)+yQ(B,\tau)+R(B,\tau).
	\end{equation}
	On the first interval $-T<\tau<T$, linearity gives the single-ratio expression
	\begin{equation}\label{eq:phase-single-ratio}
		\Phi_{x,y}(B,\tau)=\mathcal D_B\bigl[
		xA(\,\cdot+\ii(\tau-u))+yA(\,\cdot+\ii(\tau+u))
		-\ell_u(\,\cdot+\ii\tau)\bigr].
	\end{equation}
	For a fixed $(x,y)$, a \emph{critical point} is a pair $(B_*,\tau_*)$ in the appropriate space \eqref{eq:parameter-spaces} satisfying
	\begin{equation}\label{eq:critical-equations}
		\partial_B\Phi_{x,y}(B_*,\tau_*)=\partial_\tau\Phi_{x,y}(B_*,\tau_*)=0.
	\end{equation}
	
	\section{Main results}\label{sec:theorems}
	
	\begin{theorem}[Antiferroelectric regime]\label{thm:AF}
		Fix $\Delta<-1$ and $r=a/b>0$, and use parameters from \eqref{eq:AFparams}. The outer curve $\cO$ and inner curve $\cI$ recalled and defined, respectively, in Section~\ref{sec:curves}, are simple closed curves. 
		
		For every $(x,y)$ between these two curves, \eqref{eq:critical-equations} has a unique solution $(B_*,\tau_*)$ in $\mathcal P_{\rm AF}$. In this liquid region, define $h_{\rm AF}(x,y):=\Phi_{x,y}(B_*,\tau_*)$. Inside $\cI$, define $h_{\rm AF}(x,y):=(x+y)/2-T/(2\zeta)$, while outside $\cO$ define $h_{\rm AF}(x,y)$ according to affine functions \eqref{eq:frozen-planes}. These definitions extend continuously and agree on boundaries. 
		
		Then $h_{\rm AF}(x,y)$ is a limit shape. 
	\end{theorem}
	
	\begin{theorem}[Disordered regime]\label{thm:D}
		Fix $-1 < \Delta\le1/2$ and $r=a/b>0$, and use parameters from \eqref{eq:Dparams}-\eqref{eq:D-u}. The outer curve $\cO$ is recalled in Section~\ref{sec:curves}; it is a simple closed curve. 
		For every $(x,y)$ inside this curve except the center $(1/2,1/2)$, \eqref{eq:critical-equations} has a unique solution $(B_*,\tau_*)$ in $\mathcal P_{\rm D}$. Define $h_{\rm D}(x,y):=\Phi_{x,y}(B_*,\tau_*)$ in this region; additionally, define $h_{\rm D}(1/2,1/2):=u/(2 \delta)$. Outside $\cO$ define $h_{\rm D}(x,y)$ according to affine functions \eqref{eq:frozen-planes}. These definitions extend continuously and agree on boundaries. 
		
		Then $h_{\rm D}(x,y)$ is a limit shape.
	\end{theorem}
	
	\begin{remark}
		The restriction $\Delta\le1/2$ in Theorem \ref{thm:D} is purely for technical reasons: It is highly plausible that the claim holds for all values $-1 < \Delta <1$. It is plausible that the limit shape result for $\Delta=-1$ can be obtained by a relatively minor variation of our proof, but we also omit it for technical reasons. 
	\end{remark}
	
	\section{Arctic curves}\label{sec:curves}
	
	\subsection{Envelope}
	For continuously differentiable real functions $a_0,b_0,e_0$ on an interval of the real line, with $W=a_0b_0'-a_0'b_0\ne0$, define
	\begin{equation}\label{eq:envelope}
		\cE[a_0,b_0,e_0](\theta)=
		\left(\frac{e_0b_0'-e_0'b_0}{W},
		\frac{a_0e_0'-a_0'e_0}{W}\right).
	\end{equation}
	This is the unique solution $(x,y)$ of
	\begin{equation}\label{eq:envelope-equations}
		xa_0+yb_0=e_0,\qquad xa_0'+yb_0'=e_0'.
	\end{equation}
	For the analytic functions below, this is a parameterization of the envelope of the corresponding family of lines.
	
	\subsection{The outer curve in both regimes}
	This curve was predicted by Colomo-Pronko \cite{CP} and Colomo-Pronko-Zinn-Justin \cite{CPZ}. 
	Use the functions $A$ and $\ell_v$ of the corresponding regime. For $0<v<m$ define the southwest arc
	\begin{equation}\label{eq:outer-general}
		\Gamma_v(\theta)=(x_v(\theta),y_v(\theta))
		=\cE\bigl[A(\ii(\theta-v)),A(\ii(\theta+v)),\ell_v(\ii\theta)\bigr],
		\quad |\theta|<m-v,
	\end{equation}
	It is obtained in the $B\to 0$ limit. 
	
	The closed outer curve $\cO$ is the union of the following closed arcs:
	\begin{equation}\label{eq:four-arcs}
		\begin{aligned}
			\Gamma_{\rm SW}&=\{(x_u(\theta),y_u(\theta)):-T\le\theta\le T\},\\
			\Gamma_{\rm SE}&=\{(1-x_T(\theta),y_T(\theta)):-u\le\theta\le u\},\\
			\Gamma_{\rm NE}&=\{(1-x_u(\theta),1-y_u(\theta)):-T\le\theta\le T\},\\
			\Gamma_{\rm NW}&=\{(x_T(\theta),1-y_T(\theta)):-u\le\theta\le u\}.
		\end{aligned}
	\end{equation}
	
	\paragraph{Antiferroelectric coefficients.}
	For \eqref{eq:outer-general}, let us set
	\begin{equation}\label{eq:AFouter}
		A_\zeta(\ii t)=\frac{\sinh\zeta}{2\pi(\cosh\zeta-\cosh t)},\qquad
		\ell_u(\ii\theta)=\frac{T}{2\pi\zeta}
		+\frac1\pi\sum_{n\ge1}e^{-\zeta n}\frac{\sinh(Tn)}{\sinh(\zeta n)}\cosh(n\theta).
	\end{equation}
	The series converges at $\theta=\pm T$ as well. This parameterization is the Colomo-Pronko-Zinn-Justin outer-curve family \cite{CPZ}.
	
	
	\paragraph{Disordered coefficients.}
	In this case the coefficients are simpler:
	\begin{equation}\label{eq:Douter}
		\begin{aligned}
			A_\gamma(\ii t)&=\frac{\sin\gamma}{2\pi(\cos t+\cos\gamma)},\\
			\ell_u(\ii\theta)&=\frac{\sin T}{2\pi(\cos\theta-\cos T)}
			-\frac{\sin(\pi T/\delta)}{2\delta(\cos(\pi\theta/\delta)-\cos(\pi T/\delta))},\\
		\end{aligned}
	\end{equation}
	The end points of the curve are defined by corresponding limits.
	
	This curve was predicted by Colomo-Pronko \cite{CP}, and proved for $\Delta = 1/2, a/b=1$ by Aggarwal \cite{Agg20}. 
	
	\subsection{The antiferroelectric inner curve}
	Define real periodic functions by the absolutely convergent series
	\begin{equation}\label{eq:inner-fe}
		\begin{aligned}
			f(\tau)&=\frac1{2\zeta}\sum_{j\ge0}
			\frac{(-1)^j\cos((2j+1)\pi\tau/(2\zeta))}{\sinh((2j+1)\pi^2/(2\zeta))},\\
			e(\tau)&=\frac1{2\zeta}\sum_{n\ge1}
			\frac{(-1)^{n-1}\sin(n\pi u/\zeta)\cos(n\pi\tau/\zeta)}{\sinh(n\pi^2/\zeta)}.
		\end{aligned}
	\end{equation}
	Set
	\begin{equation}\label{eq:inner-PQ}
		P_1(\tau)=f(\tau-u),\quad Q_1(\tau)=f(\tau+u),\quad
		W_1(\tau)=P_1Q_1'-P_1'Q_1.
	\end{equation}
	The inner curve $\cI$ is
	\begin{equation}\label{eq:inner-curve}
		\quad
		x(\tau)=\frac12+\frac{Q_1e'-Q_1'e}{W_1},\qquad
		y(\tau)=\frac12+\frac{P_1'e-P_1e'}{W_1},\qquad \tau\in\mathbb T_\zeta.
	\end{equation}
	This curve is obtained in the limit $B \to \pi$.  
	
	\subsection{Geometry of the antiferroelectric inner curve}\label{sec:limits}
	
\begin{figure}
		\centering
		\includegraphics[width=0.97\textwidth]{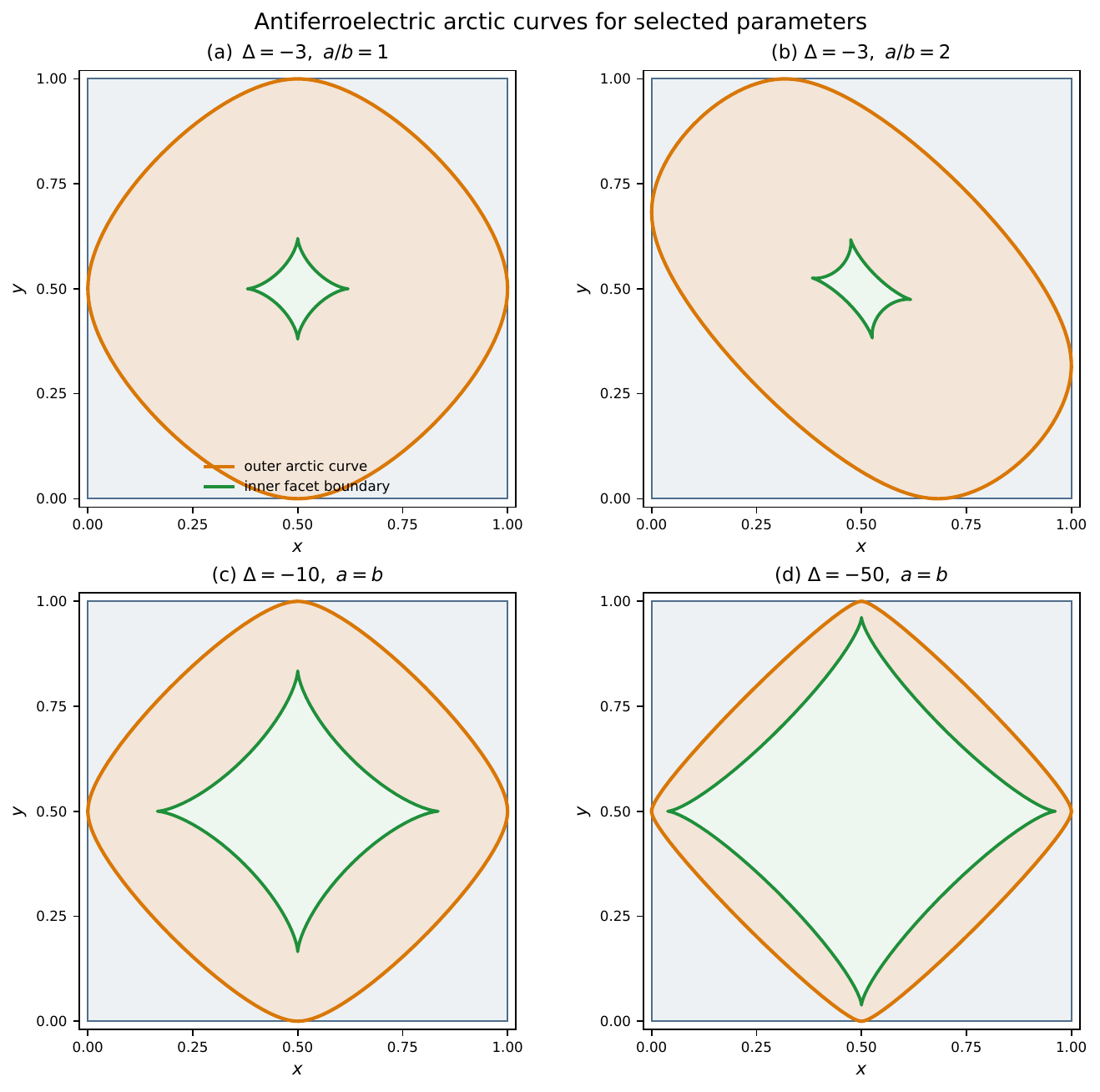}
		\caption{Examples of antiferroelectric arctic curves.}
		\label{fig:AF-requested-all}
	\end{figure}

	See Figure \ref{fig:AF-requested-all} for examples of the arctic curves in AF regime.
	
	\paragraph{Cusps in the symmetric case $a=b$.}
	In this case $u=T=\zeta/2$. The curve has four cusps. Define the complete elliptic integral
	\[
	\mathbf K(k)=\int_0^{\pi/2}(1-k^2\sin^2t)^{-1/2}\dd t,
	\qquad \tilde k=\sqrt{1-k^2},
	\]
	and choose the unique $k\in(0,1)$ such that
	\begin{equation}\label{eq:elliptic-modulus}
		\frac{\mathbf K(\tilde k)}{\mathbf K(k)}=\frac\zeta\pi.
	\end{equation}
	The four cusps are
	\begin{equation}\label{eq:cusps}
		\left( 1/2\pm \frac{\tilde k}{2},1/2 \right),\qquad \left( 1/2,1/2\pm \frac{\tilde k}{2} \right).
	\end{equation}
	
	\paragraph{Degree}
	
	The inner curve is an algebraic curve of degree at most 8. 
	
	\section{Brief History}
	
	The physical origin of the ice rule predates the six-vertex model.  In 1935 \cite{Pauling35}, Pauling used local proton-disorder constraints in ordinary ice to estimate its residual entropy.  The square-ice and six-vertex models later provided two-dimensional lattice realizations of this type of local conservation law.  Their exact solution was initiated by Lieb and Sutherland in the 1960s \cite{Lieb67,Suth67}.
	
	The choice of \emph{domain-wall boundary conditions} (DWBC) was introduced by Korepin \cite{Korepin82}.  The exact determinant formula for the corresponding finite-volume partition function was found by Izergin \cite{Izergin87}.  The Izergin-Korepin determinant became the essential algebraic tool for the DWBC model and later played an important role in both asymptotic analysis and enumerative combinatorics.
	
	A second classical source of interest is the relation with alternating sign matrices.  The connection was already implicit in the work of Robbins and Rumsey \cite{RR86} and in the study of Aztec-diamond domino tilings by Elkies, Kuperberg, Larsen, and Propp \cite{EKLP92}.  It became central after Kuperberg's proof of the alternating-sign-matrix conjecture via the six-vertex model with DWBC \cite{Kuper96}; see also Zeilberger's proof \cite{Zeil96} and Zinn-Justin's matrix-model and determinant reformulations \cite{ZJ00}.
	
	The point $\Delta=0$, equivalently $c^2=a^2+b^2$, is the \emph{free-fermion point} and has a particularly rich additional structure.  The corresponding determinantal model is closely related to the dimer model, whose exact solution goes back to Temperley-Fisher and Kasteleyn \cite{TF61,Kast61}.  For DWBC, the free-fermion six-vertex model can be mapped to weighted domino tilings, see Ferrari-Spohn \cite{FS06}.  The pioneering works on dimers include the exact Aztec-diamond enumeration of Elkies-Kuperberg-Larsen-Propp \cite{EKLP92}, the arctic-circle theorem of Jockusch--Propp--Shor \cite{JPS98}, the local-statistics and limit shape results of Cohn-Elkies-Propp \cite{CEP96}, Cohn-Kenyon-Propp \cite{CKP01}, Kenyon-Okounkov \cite{KO07}, and Johansson's Airy-process fluctuations of the arctic boundary \cite{Joh05}.
	
	The numerical evidence for a nontrivial macroscopic separation of regions in the DWBC model was observed, in particular, by Allison and Reshetikhin \cite{AR05}.  Large-$N$ asymptotics of the partition function were obtained by Bleher and Fokin in the disordered phase \cite{BF06}, and subsequently by Bleher and Liechty on the ferroelectric side, on the critical line between ferroelectric and disordered phases, and in the antiferroelectric phase \cite{BL08,BL09,BL10}.
	
	The first explicit formulas for the outer frozen boundary in the DWBC geometry were predicted by Colomo and Pronko in the disordered regime \cite{ColomoPronkoLimit10,CP} and by Colomo, Pronko, and Zinn-Justin in the antiferroelectric regime \cite{CPZ}.  Later, Colomo and Sportiello developed the tangent method, which gave a conceptual derivation of these curves for a wide class of domains \cite{CS16}.  A complementary variational perspective on fixed-boundary six-vertex models was proposed by Reshetikhin and Palamarchuk \cite{RP10}.
	
	Variational principles and large-deviation statements for random height functions have a rich history in random-surface models. Particularly important works include Cohn-Kenyon-Propp \cite{CKP01}, Sheffield \cite{Shef05}, Lammers-Tassy \cite{LT} and Kenyon-Prause \cite{KenyonPrause22,KenyonPrause24}. 
	
	At the ice point $\Delta=1/2$, $r=1$, Aggarwal gave a rigorous proof of the arctic-curve description predicted earlier by Colomo and Pronko; see \cite{Agg20}.  Beyond the law of large numbers, Ayyer, Chhita, and Johansson proved a first non-free-fermionic edge-fluctuation theorem for the domain-wall model \cite{ACJ23}.  More recently, Gorin and Liechty studied boundary statistics for general DWBC weights, proving GUE-corners asymptotics for $\Delta<1$ and finite-order stochastic-six-vertex limits for $\Delta>1$ \cite{GL25}. 
	
	\subsection*{Comments}
	
	The proof of the result was assisted by ChatGPT.
	
	The proof follows the general strategy from the previous work in this series (announced in \cite{BO}). It uses ideas and results, among others, from Cohn-Kenyon-Propp \cite{CKP01}, Sheffield \cite{Shef05}, Lammers-Tassy \cite{LT}, Alessandrini \cite{Alessandrini87}, Duminil-Copin-Kozlowski-Krachun-Manolescu-Tikhonovskaia \cite{DCK}, Kitanine-Maillet-Terras \cite{KMT99}, and Kenyon-Prause \cite{KenyonPrause22,KenyonPrause24}. The conceptual description of Kenyon-Prause \cite{KenyonPrause22,KenyonPrause24} is especially relevant for the results described above. 
	
	Additional technical difficulties that arise in this setup compared to the square ice include the absence of previously rigorously proven arctic curve results and the significant influence of a new antiferroelectric region on variational arguments.

	\subsection*{Acknowledgements}
	
	A.~Bufetov was partially supported by the European Research Council (ERC), Grant Agreement No. 101041499.

	(Alexey Bufetov) Institute of Mathematics, Leipzig University, Germany
	
	Email: alexey.bufetov@gmail.com
	
\end{document}